%Paper: dg-ga/9407015
%From: mmurray@maths.adelaide.edu.au (Michael Murray)
%Date: Mon, 25 Jul 1994 16:13:36 +0930

\magnification =1200
\footline={\hss\tenrm\folio\hss}
\normalbaselineskip = 1.2\normalbaselineskip
\normalbaselines
\parindent=20pt
\overfullrule=0pt

\font\titlefont=cmbx10 scaled \magstep4
\font\authorfont=cmbx10 scaled \magstep2

\def\C{{\bf C}}
\def\Z{{\bf Z}}
\def\R{{\bf R}}

\def\a{\alpha}
\def\b{\beta}
\def\c{\gamma}

\def\ref#1{\par\hangindent=1.0cm\hangafter1{\noindent#1}}
\def\tilde{\widetilde}
\def\Cstar{{\bf C}^{\times}}

\def\ev{\mathop{\rm ev}\nolimits}
\def\wzw{\mathop{\rm wzw}\nolimits}

\def\hol{\mathop{\rm hol }\nolimits}
\def\and{\hbox{\quad and \quad}}

%
% Title page
%

\nopagenumbers

\null
\rightline{dg-ga/9407015}
\vskip 1.5 cm
\centerline{\titlefont Bundle gerbes}
\vskip .5 cm
\vskip 1 cm
\centerline{\authorfont M.K. Murray$^1$}
\vskip 1 cm
%\centerline{\authorfont }
\vskip 1.5cm

\vfill

\centerline{\bf 25 July 1994}

\vfill

\noindent 1. Department of Pure Mathematics, The University of Adelaide,
Adelaide SA 5005, Australia.  {\it mmurray@maths.adelaide.edu.au}

\vfill\eject
\nopagenumbers

\vskip 1cm

\noindent{\bf Abstract.} Just as $\Cstar$ principal
bundles provide a geometric realisation of two-dimensional
integral cohomology; gerbes  or sheaves of groupoids,
provide a geometric realisation
of three dimensional integral cohomology through their Dixmier-Douady
class. I consider an
alternative, related,  geometric realisation of three
dimensional cohomology called a bundle gerbe. Every bundle
gerbe gives rise to a gerbe and most of the well-known examples
examples of gerbes are bundle gerbes. I discuss the properties of
bundle gerbes, in particular bundle gerbe
connections and curvature and their associated Dixmier-Douady class.

\vskip 1cm
\noindent{\bf Mathematics Subject Classification:}
55R65, 18B40.

\vfill\eject

\pageno=1
\footline={\hss\tenrm\folio\hss}

\bigskip
\noindent{\bf  1. Introduction \hfill}

In [1] Brylinski describes Giraud's theory of gerbes. Loosely speaking
a gerbe over a manifold $M$ is a sheaf of groupoids over $M$. Gerbes,
via their Dixmier-Douady class,
provide a geometric realisation of the elements of $H^3(M, \Z)$ analogous
to the way that line bundles provide, via their Chern class,
a geometric realisation of the elements of $H^2(M, \Z)$.

I want to introduce here another sort of object,
which is not a sheaf, and which also gives rise to elements of
$H^3(M, \Z)$. For want of a better name I have called these
objects bundle gerbes. A bundle gerbe over $M$ is a pair
consisting of a fibration $Y \to M$ and a principal $\Cstar$  bundle
$P$ over the  fibre product $Y^{[2]}$. The bundle $P$ is required to
have a product,  that is a  $\Cstar$  bundle morphism which on fibres
is of the form:
$$
P_{(x,y)} \otimes P_{(y,z)} \to P_{(x,z)},
\eqno(1.1)
$$
for any $x,y,z$ in the same fibre of $Y \to M$.
{}From a bundle gerbe it is possible to  construct a
presheaf of groupoids and hence a gerbe. However not all gerbes
arise in this way.

By considering a  connection on  $P$ compatible with the product  (1.1)
it is possible to  construct a closed, integral three-form
 on $M$ analogous to the
curvature of a line bundle.  This  three form is a representative
for the image in $H^3(M, \R)$ of the
Dixmier-Douady class of the corresponding gerbe. Moreover
every integral three  class is represented by a closed, integral three
form arising from a bundle gerbe $P \to Y^{[2]}$ where $Y \to M$
is the path fibration.  This construction is analogous to
the construction of a line bundle with given curvature two-form.
It follows from the results in [1] that every isomorphism class
of gerbes contains a bundle gerbe.
If $Q \to Y$ is a $\Cstar$ bundle there is  a bundle gerbe whose
fibre at $(x,y)$ is $Aut_{\Cstar}(Q_x, Q_y)$.
I shall call such a bundle gerbe  trivial.
The geometric interpretation of the  Dixmier-Douady class of a bundle
gerbe is that it is the obstruction to the bundle gerbe being trivial.

Many interesting gerbes are bundle gerbes.
For instance if $\C^\times \to \hat G \to G$ is a
central extension of groups
it is well known that the obstruction to lifting a principal $G$
bundle over $M$ to a principal $\hat G$ bundle  is a class in
$H^3(M, \Z)$ and Brylinski discusses the gerbe defined by such
a principal bundle whose Dixmier-Douady class is this
obstruction. This principal bundle also gives rise to
a bundle gerbe in a natural way and  this bundle
gerbe is trivial precisely when the original bundle lifts
to $\hat G$.

Having outlined the virtues of bundle gerbes  I should
mention two of their deficiencies.  The first
is that there are many bundle gerbes which are isomorphic as
gerbes but not as bundle gerbes. As a consequence the theorem that
isomorphism  classes of gerbes are in bijective correspondence with
$H^3(M, \Z)$ does not  hold for bundle gerbes.   The second is that
if a bundle gerbe $P \to Y^{[2]}$ is  non-trivial  then the fibres of
$Y \to M$ must  be infinite dimensional.

In outline the paper is as follows.  Section 2 reviews the
properties of $\Cstar$ groupoids and their relationship with
line bundles. This is in preparation for Section 3 which
introduces bundle gerbes. The motivating example of the lifting of
a principal bundle for a central extension is discussed in Section 4
and this leads to the introduction of the Dixmier-Douady class of a
bundle gerbe in Section 5.  A de Rham representative for this class is
provided by the theory of bundle gerbe connections and curvature introduced
in Sections 6, 7 and 8. The relationship with gerbes
 is discussed in Section 9 as in the relationship of the bundle
gerbe connection and curvature with Brylinski's connective
structure and curving. In Section 10 I show how to construct
a bundle gerbe  with given Dixmier-Douady class.
Section 11 considers the Deligne cohomology class defined
by a bundle gerbe with connection and in Section 12  the
holonomy of a bundle gerbe connection over a  two-sphere is defined.
Finally in Section 13 I explain why the fibering $Y \to M$ has to
have infinite dimensional fibres if the bundle gerbe is non-trivial.

I will assume, when talking
about gerbes, some familiarity with Brylinski's book [1]. The material
on bundle gerbes however is intended to be self-contained.

\bigskip
\noindent{\bf   2. $\C^\times$ groupoids \hfill}

Denote by $\Cstar$ the group of non-zero complex numbers. If $P$ and
$Q$ are principal $\Cstar$ bundles  over a manifold $Z$ then it is
possible to define a new principal $\Cstar$ bundle $P\otimes Q$ over
$Z \times Z = Z^2$. This is called the contracted product [1]. If $P$
and $Q$ are the frame bundles of  line bundle $L$ and $J$ respectively
then $P\otimes Q$ is the frame bundle of the line bundle $L\otimes
J$. To define $P\otimes Q$ take  $P\times Q$ which is a  principal
$\Cstar \times \Cstar$ bundle over $Z^2$ and  quotient by the
`anti-diagonal' copy of $\Cstar$ inside $\Cstar\times \Cstar$, that
is,  the subgroup of all pairs $(z, z^{-1})$. What makes this
construction possible, of course, is the fact that  $\Cstar$ is
abelian. Consider now a manifold $X$ and inside $X^2 \times X^2$
define
$$
X^2 \circ X^2 = \{ ((x, y), (y,z)) \mid x, y, z \in X\}.
\eqno(2.1)
$$
If $P$ is a principal $\Cstar$ bundle over $X^2$ define $P\circ P$
to be the restriction of $P\otimes P$ to  $X^2 \circ X^2$.

Recall [2] that a groupoid is  a category with every morphism
invertible. Let us  consider an equivalent definition
that we will see below is easily generalised to define
bundle gerbes.  Define a  $\Cstar$ groupoid  to be a
principal $\Cstar$ bundle $P$ over $X^2$
with a product, that is, a $\Cstar$ bundle morphism $P\circ P \to P$,
$(p,q) \mapsto pq$, covering the map $((x,y),(y,z)) \mapsto (x,z)$.
The product is required to be  associative, that is $(pq)r = p(qr)$ whenever
these products are defined.

A $\Cstar$ groupoid  actually has two other important algebraic
structures, an identity
and inverse  which could have been included in the definition but in fact
are a consequence of it.   The identity is a section $e$
of the bundle $P$ over the diagonal in $X^2$ which satisfies $pe = ep = p $.
To define it
note first that if $p \in P_{(x,y)}$  and $q \in P_{(y,y)}$
then $pq \in P_{(x,y)}$. Hence there is some $z \in \Cstar$
such that $pq = pz$. Define $e = qz^{-1}$ so that
$pe = p$. Because the product is a bundle automorphism
 $(pw)e = (pe)w = pw$ for all $w \in \Cstar $ and
hence $qe = q$ for every $q$ in $P$.  To show that
 $ep=p$ we use associativity. Clearly
$ep =  pz$ for some $z \in \C^\times$ and  considering
$(pe)p = p(ep)$ it follows that  $z = 1$.
To define the inverse   notice that  the equation
$pq = e$ can be solved by by acting by $\Cstar$. Then we have
that $qp = ez$ for some $z$ and
 using associativity  in the form of $p(qp) = (pq)p$ it follows that
$pez = ep = pe$ and hence $z = 1$. The inverse will be denoted
by   $p \mapsto p^{-1}$.
To understand the global structure of the inverse notice that
it is possible to construct a bundle $P^*$ over $X^2$ by defining it
to be the  same set as $P$ but changing the $\Cstar$ action to
$pz = pz^{-1}$. If $P$ is the frame bundle of a line bundle
$L$ then $P^*$ is the frame bundle of $L^*$.
Because $\Cstar$ is abelian this is still a
right action. The inversion then defines a map $P \to P^*$
covering the map $X^2 \to X^2$ defined by $(x,y) \mapsto (y,x)$.
The identity and the inverse behave as one would  expect with respect
to the product.

Given a $\Cstar$ groupoid  we can recover the
definition in terms of   categories [2] by
taking $X$ as the set of objects and $P_{(x,y)}$ as the morphisms
from $x$ to $y$.

A simple example of a $\Cstar$ groupoid is constructed by taking a
$\Cstar$ bundle $Q$ on $X$ and defining $P_{(x,y)} = Aut_{\Cstar}(Q_x,
Q_y)$ where $Q_x$ is the fibre of $Q$ over $x$ and the subscript
$\Cstar$ indicates that these automorphisms commute with the $\Cstar$
action. An alternative way to define theis is to use the two
projections $\pi_1$ and $\pi_2$ on the first and second factors of
$X^2 = X\times X$ and define $P = \pi_1^{-1}Q^*\otimes \pi_2^{-1} Q$.
The composition is, of course, composition of automorphisms. What
makes $\Cstar$ groupoids uninteresting  is that every $\Cstar$
groupoid arises in this way!  To see this let $P$ be a $\Cstar$
groupoid and  pick a basepoint $x$ in $X$. Then define a $\Cstar$
bundle $Q$ on $X$ by $Q_y = P_{(x,y)}$, that is $Q$ is the pull-back
of $P$ under the map $y \mapsto (x,y)$. Then the composition and
inversion can be used to define a $\Cstar$ bundle isomorphism
$$
Q_y^* \otimes Q_z = P_{(x,y)}^* \otimes P_{(x,z)} \to P_{(y,z)}
\eqno(2.2)
$$
by $(p,q) \mapsto p^{-1}q$.  It is easy to see that this is a $\Cstar$
bundle isomorphism and that, moreover, it preserves the composition.
Hence it is  an isomorphism of $\Cstar$ groupoids.

Although we have just seen that the theory of $\Cstar$
groupoids is nothing more than the theory  of $\Cstar$ bundles
over pointed sets it is useful to develop the theory further
as the next section of bundle gerbes is then a straightforward
generalisation.

A connection $\nabla$ on the bundle $P \to  X^2$ gives rise to a connection
on $P\otimes P$ and hence on $P\circ P$. Call it a
groupoid connection if it is mapped by the product to itself again.
Such connections exist  because we can identify $P$ with
$ \pi_1^{-1}Q^*\otimes \pi_2^{-1} Q$ for some $Q \to X$
and pick a connection on $Q$ and pull it back to a connection on
$P$.  The curvature $F_\nabla$ of a groupoid connection on $P$ constructed
in this way has the form
$$
F_\nabla = \pi_2^*(f) - \pi_1^*(f)
\eqno(2.3)
$$
where $f$ is the curvature of the original connection on $Q \to X$.
It is also possible to show (see Section 8 below) that if $F_\nabla$
is the curvature of any groupoid connection then there is
unique two-form $f$ on $X$ which  satisfies (2.3). Call this
two-form $f$ the groupoid curvature.

It is well known that given an integral closed two-form $f/(2\pi i)$
over a $1$-connected space $X$ it is possible to explicitely construct a
$\Cstar$ bundle over $X$ with a connection whose curvature is $f$.  It
will be useful later to know  how to repeat this construction for
$\Cstar$ groupoids. In fact this construction is a little more natural
in the groupoid context as, unlike the construction of the $\Cstar$
bundle over $X$,  it does not require the choice of a base point.

Assume then that  $ X$ is a $1$ connected manifold and that $f/(2\pi
i)$ is a closed, integral  two-form on $X$. We shall construct a
$\Cstar$ groupoid on $X$  with a groupoid connection whose groupoid
curvature is $f$.  Let $PX$ be the   space  of all piecewise smooth
paths in $X$ and let $PX_{(x,y)}$ be the set of all paths beginning at
$x$ and ending at $y$. We need piecewise smooth paths as we will want
to compose them to define the groupoid product. Define an equivalence
relation on $PX_{(x,y)} \times \Cstar$ by saying that $(\gamma, z) $
is equivalent to $(\tilde \gamma, \tilde z)$ if
$$
z = \exp(\int_D f) \tilde z
\eqno(2.4)
$$
where $D$ is a map of the disk into $X$ with boundary the union of
$\gamma $ and $\tilde \gamma$ and oriented by $\gamma$. Let
$P_{(x,y)}$ denote the quotient space. The $\Cstar$ action is
just that induced by $(\gamma, z)w = (\gamma, zw)$. If
$\gamma(1) = \tilde \gamma(0)$ then  define a new piecewise smooth
path $\gamma\star\tilde\gamma$ by running along $\gamma$ at twice the
speed for $ t\in [0, 1/2]$ and then running along $\tilde\gamma$ at
twice the speed for $t \in [1/2, 1]$. A product on $PX \times \C$ is
then defined by $(\gamma, z)\star (\tilde \gamma, \tilde z) =
(\gamma\star\tilde\gamma, z\tilde z)$ and it is straightforward to
check that this descends to a product on $P$ making it a $\Cstar$
groupoid.

Now we construct a $\Cstar$ groupoid connection. Notice that $P$ is
the quotient of a fibering
$$
PX \times \C \to P
\eqno(2.5)
$$
where the fibres are defined by the equivalence relation (2.4).
The tangent space to $PX \times \C$ at $(\gamma, z)$
is the space of all pairs
$(\xi, \alpha)$ where $\xi$ is a vector field along $\gamma$ and
$\alpha$ is a complex number. A vector field
along $\gamma$  means a continous vector field which is smooth
when $\gamma$ is smooth. The subspace of the tangent
space that is tangent to the fibering (2.5) is
$$
K_{(\gamma, z)} = \{ (\xi, \alpha) \mid \alpha =
- \int_0^1 f(\gamma', \xi)dt,  \xi(0) = 0 = \xi(1) \}
\eqno(2.6)
$$
where $\gamma'$ is the tangent vector field along $\gamma$.
Consider the map
$$
\ev \colon PX \times [0,1] \to X
\eqno(2.7)
$$
which maps $(\gamma, t) \mapsto \gamma(t)$.  A  one-form $\hat A$
on $PX \times [0,1]$ is now defined by pulling back $f$ and integrating
it over the  $[0,1]$ direction and then letting
$$
\hat A = dt + \int^1_0 \ev^*(f).
\eqno(2.8)
$$
Notice that
$$
d\hat A = \ev_1^*(f) - \ev_0^*(f)
\eqno(2.9)
$$
where $\ev_t(\gamma) = \gamma(t)$.  The forms $\hat A$ and $d\hat A$
both annihilate vectors in the space (2.6) and hence $\hat A$ descends
to a one-form $A$ on $P$ which defines a  connection.

It remains to check that this connection is a groupoid connection. The
product map defines a sum on tangent vectors. If $(\xi, \alpha)$, and
$(\tilde\xi, \tilde\alpha)$ are tangent to $(\gamma, z)$ and
$(\tilde\gamma, \tilde z)$ then $(\xi\star\tilde\xi , \alpha+
\tilde\alpha)$ is tangent at $(\gamma, z)\star(\tilde\gamma, \tilde
z)$ where $\xi\star\tilde\xi$ is the obvious vector field along
$\gamma\star\tilde\gamma$. The essential point in the proof that $A$
preserves the product is that
$$
\int_0^1 f(\gamma', \xi) dt + \int_0^1 f(\tilde\gamma', \tilde\xi) dt
= \int_0^1 f({(\gamma\star\gamma)'}, \xi\star\tilde\xi) dt.
\eqno(2.10)
$$
It follows from (2.9) that
$$
dA = \pi_2^*(f) - \pi^*_1(f)
\eqno(2.11)
$$
and hence the curvature of this groupoid connection is $f$.

In the next section we are concerned with  bundle gerbes
which are  fibrations whose fibres are $\Cstar$ groupoids.
Then it may  not possible to choose basepoints continously
and the constructions above become more interesting.

\bigskip
\noindent{\bf  3. $\Cstar$ bundle gerbes \hfill}

Consider a fibration  $\pi \colon Y \to M$. Define the fibre product
$Y^{[2]}$ in the usual way, that is the subset of pairs $(y, y')$ in
$Y \times Y$ such that $\pi(y) = \pi(y')$. Notice that the
diagonal is inside $Y^{[2]}$ and that  the map that transposes elements
of $Y^2 = Y \times Y$ fixes $Y^ {[2]}$. Denote by $\pi_i$ the restriction
of the projection maps on
$Y^2$ to $Y^{[2]}$. Denote by $Y^{[2]}\circ Y^{[2]}$ the intersection
of $Y^{[2]}\times Y^{[2]}$ with $Y^2\circ Y^2$. If $P$ and $Q$ are $\Cstar$
bundles over $Y^{[2]}$ denote by $P\circ Q$ the restriction to
$Y^{[2]}\circ Y^{[2]}$ of the bundle $\pi_1^{-1}(P) \otimes
\pi_2^{-1}(Q)$ over $Y^{[2]} \times Y^{[2]}$.

A bundle gerbe over $M$ is defined to be a choice of a fibration $Y
\to M$ and a $\Cstar$ bundle $P \to Y^{[2]}$ with a product, that is,
a $\Cstar$ bundle isomorphism $P\circ P \to P$ covering $((y_1,
y_2),(y_2, y_3)) \mapsto (y_1, y_3)$. The product is required to be
associative whenever triple products are defined. Just as for $\Cstar$
groupoids  a bundle gerbe has an inverse and an identity denoted by
the same symbols. Occasionally we shall denote a bundle gerbe as a
triple $(P, Y, M)$.

\medskip

\noindent{\bf Example: 3.1 } Let $Q \to Y$ be a principal $\Cstar$ bundle.
Define $P_{(x, y)} = Aut_{\Cstar}(Q_x, Q_y) = Q_x^* \otimes Q_y$
Then this defines a  bundle gerbe  called the trivial bundle gerbe.
We also have $P = Aut(\pi_1^{-1}Q, \pi_2^{-1}{Q}) = \pi_1^{-1}Q^*
\otimes \pi_2^{-1}Q$.

\medskip

\noindent{\bf Example: 3.2 } An example of interest in [1] is to start
with a fibration $Y \to M$ with $1$ connected fibres and a two-form on $Y$
whose restriction to each fibre is closed and integral. Then we
can apply the construction in Section 2 fibre by fibre to
define a bundle gerbe $P \to Y^{[2]}$.

\medskip

A morphism of bundle gerbes $(P, Y, M)$ and $(Q , X , N)$ is a triple
of maps $(\a, \b, \c)$. The map $\beta \colon Y \to X$ is required to
be a morphism of the fibrations $Y \to M$ and $X \to N$ covering
$\gamma \colon M \to N$. It therefore induces a morphism $\beta^{[2]}$
 of the fibrations $Y^{[2]} \to M$ and $X^{[2]} \to N$. The  map
$\alpha$ is required to be a morphism of $\Cstar$ bundles covering
$\beta^{[2]}$  which commutes with the product and hence also with
the identity and inverse.
A morphism of bundle gerbes over $M$ is a morphism of bundle
gerbes for which $M = N$ and $\gamma$ is the identity on $M$.

Various constructions are possible with bundle gerbes.
We can define a pull-back and product as follows.
If $(Q, X, N)$ is a bundle gerbe and $f \colon M \to N$
is a map then we can pull back the fibration $X\to N$ to
obtain a fibration $f^{-1}(X) \to M$ and a morphism of
fibrations $f^{-1} : f^{-1}(X) \to X$ covering $f$.
This induces a morphism $(f^{-1}(X))^{[2]} \to X^{[2]}$ and
hence we can use this to pull back the $\Cstar$ bundle
$Q$ to a  $\Cstar$ bundle $f^{-1}(Q)$ say on $f^{-1}(X)$.
It is easy to check that $(f^{-1}(Q), f^{-1}(X), M)$
is a bundle gerbe, the pull-back by $f$ of the gerbe
$(Q, X, N)$. If $(P, Y, M)$ and $(Q, X, M)$
are bundle gerbes over $M$ then we can form a fibre
product $Y\times_M X  \to M$ and then form a $\Cstar$ bundle
$P \otimes Q$ over $(Y\times_M X )^{[2]}$  which is the
product of the bundle gerbes $(P, Y, M)$ and $(Q, X, M)$.

Notice that, unlike the case of $\Cstar$ groupoids, it is not
clear that every bundle gerbe is trivial.  The proof in
Section 2 that a groupoid is the same as a bundle over a pointed set
can only be applied fibre by fibre if $Y \to M$ has a
section.  We shall see in Section 5 that there is associated to
a bundle gerbe a class in $H^3(M, \Z)$, its Dixmier-Douady class,
which is precisely the obstruction to the bundle gerbe being
trivial.

\bigskip
\noindent{\bf 4. Central extensions \hfill}

A motivating  example of a bundle gerbe is the
bundle gerbe arising from a central extention of
groups.  Let
$$
0 \to \Cstar {\buildrel\iota \over \to } \  \hat G
{\buildrel p\over \to }\  G \to 0
\eqno(4.1)
$$
be a central extension of groups and
 $Y \to M$  a principal $G$ bundle. Then it may happen that there is
a principal $\hat G$ bundle $\hat Y$ and a bundle map
$\hat Y \to Y$ commuting with the homomorphism $\hat G \to G$.
In such a case  $Y$ is said to lift to a $\hat G$ bundle.
One way of answering the question of when $Y$ lifts to a $\hat G$ bundle
is  to present $Y$ with transition functions $g_{\a\b}$ relative to a
cover $\{ U_\a\}$ of $M$. If the cover is sufficiently nice we can
lift the $g_{\a\b}$ to maps $\hat g_{\a\b}$ taking values in
$\hat G$ and such that $p(\hat g_{\a\b}) = g_{\a\b}$.
These are candidate transition functions for a lifted bundle $\hat Y$.
However they may not satisfy the cocycle condition
$\hat g_{\a\b} \hat g_{\b\c} \hat g_{\b\c} = 1$
and indeed there is a $\Cstar$ valued map
$$
e_{\a\b\c} \colon U_\a \cap U_\b \cap U_\c \to \Cstar
\eqno(4.2)
$$
defined by $\iota(e_{\a\b\c}) = \hat g_{\a\b} \hat g_{\b\c} \hat g_{\c\a}$.
Because (4.1) is a central extension it follows  that  $e_{\a\b\c}$
is a co-cycle and hence defines a class in $H^2(M, \Cstar)$.
It is well-known that the coboundary map in the long exact
sequence of cohomology induced by (4.1) defines an isomorphism
$$
H^2(M, \Cstar) \simeq H^3(M , \Z).
\eqno(4.3)
$$
The image of $e_{\a\b\c}$ under this coboundary is the
class in $H^3(M, \Z)$ which is the obstruction to
$Y \to M$ lifting to $\hat G$.

Another way to see when $Y \to M$  lifts to $\hat G$ is to
construct a bundle gerbe which is trivial precisely when a lift
is possible. To do this
define $ P \to Y^{[2]}$ by
$$
\hat P_{(x, y)} = \{ h \in \hat H \mid x p(h)  = y\}.
\eqno(4.4)
$$
Assume now that $Y$ has a lift to a principal $\hat G$
bundle $\hat Y$ over $M$ so that there is a projection
$q \colon \hat Y \to Y$ commuting with $p$ in the appropriate way.
Then $\hat Y\to Y$ is a $\C^\times$ bundle over $Y$. Indeed
let  $g \in P_{(x, y)}$;
then it defines a map $\hat Y_x \to \hat Y_y$
which, by centrality, commutes
with the $\C^\times$ action. This defines an isomorphism
$$
P_{(x, y)} \simeq Aut_{\Cstar}(\hat Y_x, \hat Y_y)
\eqno(4.5)
$$
so that $P\simeq (\pi^*_1{\hat Y})^* \otimes \pi^*_2{\hat Y}$.
On the other hand  if the bundle gerbe
$P$ is trivial, say isomorphic to $Aut(\hat Y, \hat Y)$ for
some $\C^\times$ bundle $\hat Y \to Y$ it is possible to define an action
of $\hat G$ on $\hat Y$ and make it a lift of $Y$.  To do this start
with $g$ in $\hat G$ and the fibre $\hat Y_y$. Then define $x$ by
 $x  p(g)= y$.
Then $g \in P_{(x, y)} = Aut_{\C^\times}(\hat Y_x, \hat Y_x)$
so apply the corresponding automorphism to any element in $\hat Y_x$
to define the action of $g$.  It can be checked that this
defines a lift of $Y$.

This proves that $Y$ lifts to $\hat G$ if and only if the gerbe $P$
is trivial. In other words the bundle gerbe $P$ is trivial when the
three class defined by $Y$ is zero. We shall see in the
next section that this examples generalises to
all bundle gerbes.  Every bundle gerbe defines a three class,
its  Dixmier-Douady class, which is the obstruction to
it being trivial.

\bigskip
\noindent{\bf  5. The Dixmier-Douady class of a bundle gerbe.\hfill}

Let $P \to Y^{[2]}$ be a bundle gerb. Choose a  cover $\{U_\a\}$ of $M$
such that over
each $U_\a$ there is a section $s_\a$  of $Y$. Then
on the overlap $U_\a \cap U_\b$ we have a map
$$
(s_\a, s_\b) \colon U_\a \cap U_\b \to Y^{[2]}
\eqno(5.1)
$$
defined by $(s_\a, s_\b)(x) = (s_\a(x), s_\b(x))$.
Let $P_{\a\b}$ be the pull-back of $P$ via this map.
Notice that the product gives an isomorphism
 $P_{\a\b}\otimes P_{\b\c} \simeq P_{\a\c}$.
Choose sections $\sigma_{\a\b}$
of each $P_{\a\b}$. Then using the
product  we define a $\C^\times$ valued function
$g_{\a\b\c}$ defined by
$$
\sigma_{\a\b} \sigma_{\b\c} =
\sigma_{\a\c}g_{\a\b\c}.
\eqno(5.2)
$$
It is easy to check that $g$ defines a class in $H^2(M, \C^\times)$.
The  image of $g$ under this the isomorphism (4.3) is
called the  Dixmier-Douady class of the bundle gerbe $P$.

I claim that the Dixmier-Douady class is precisely the
obstruction to a gerbe being trivial.
To see this note  first that if $P$ trivial, say  $P = \pi_1^{-1}Q^*
\otimes \pi_2^{-2}Q$ for some bundle $Q$ on $Y$
we can define $Q_\a = s_\a^*(Q)$ and we then
have canonical isomorphisms
$$
P_{\a\b} = Q_\a^* \otimes Q_\b
\eqno(5.3)
$$
commuting with products. Hence in the construction of the cocycle
in Section 4, if we choose $\delta_\a$ to be a section
of $Q_\a$ and define $\sigma_{\a\b} =
(\delta_\a)^{-1}\otimes \delta_b$ we obtain a trivial cocycle $g$.

If on  the other hand  $g$ is trivial, say
$$
g_{\a\b\c} = \rho_{\a\b} \rho_{\b\c} \rho_{\c\a}
\eqno(5.5)
$$
where $\rho$ is $\C^\times$ valued, then we can divide $\sigma_{\a\b}$
in equation (5.2) by $\rho_{\a\b}$ and hence without loss of
generality assume that $g$ is identically one. Let $Y_\a =
\pi^{-1}(U_\a)$. Define a principal bundle $Q_\a$ over $Y_\a$ by
defining its fibre at $y$ to be
$$
(Q_\a)_y = P_{(y, s_\a(\pi(y)))}.
\eqno(5.6)
$$
The $\sigma_{\a\b}$ are elements of
$$
 P_{(s_\a(\pi(y)), s_\b(\pi(y)))} =
 {P^*_{s_\a(\pi(y)),y)}}\otimes  P_{s_\a(\pi(y)), y)}
= (Q_\a)_y^* \otimes (Q_\b)_y.
\eqno(5.7)
$$
The $\sigma_{\a\b}$ therefore define  automorphisms
between $Q_\a$ and $Q_\b$ over $Y_\a \cap Y_\b$.
The standard clutching construction  now defines a bundle $Q$ over all of
$Y$ and it is straightforward to check that this trivialises
the gerbe $P$ over $Y$.

\medskip

\noindent{\bf Example 5.1:} If the fibration $Y \to M$ admits a
global section  then we can
smooothly pick a base point in each fibre. The results of
Section 2 can then be applied to show that the bundle gerbe is
trivial. It is a trivial consequence of the definition in
Section 4 that  the Dixmier-Douady class is  zero.

The Dixmier-Douady class behaves naturally with respect to the
operations on bundle gerbes defined in section 3. The
Dixmier-Douady class of a pull-back bundle gerbe is the
pull-back of the Dixmier-Douady class and the Dixmier-Douady
class of a product of two bundle gerbes is the product (sum) of the
Dixmier-Douady classes of the individual bundle gerbes.

\bigskip
\noindent{\bf  6. Connections on bundle gerbes.\hfill}

Consider a connection $\nabla$ on a bundle gerbe $P \to Y^{[2]}$.
Then it induces a connection $\nabla\circ\nabla$ on the bundle
$P\circ P \to Y^{[2]}\circ Y^{[2]} $.
The connection $\nabla$ is said to be a  bundle gerbe connection if the
image of $\nabla\circ\nabla$ under the product  map is $\nabla$.

It is not clear that bundle gerbe connections exist. To see that they
do note that if the bundle gerbe is trivial i.e. $P = Aut(\pi_1^{-1}Q,
\pi_2^{-1}Q)$ for  some principal bundle $Q$ on $Y$ then a connection
$\nabla$ on $Q$ defines a bundle gerbe connection $\nabla^* \otimes
\nabla$  on $P$. Now  choose an open
cover $\{ U_\a\}$ of $M$ over which the
fibration $Y \to M$ is trivial.  Denote by $Y_\a$ the  open
subset of $Y$ which is the
pre-image under the projection of $U_\a$. Then the bundle gerbe can be
trivalised over $Y_\a$ and hence admits a bundle gerbe connection.
Choose a partition of unity for the open cover $\{ U_\a\}$ of $M$.  This
pulls-back to $Y^{[2]}$ to give a partition of unity for the open cover
$\{Y_\a^{[2]}\}$ and  can be used to patch together the bundle gerbe
connections on the various open sets to give a bundle gerbe connection
on $P$. It will follow  from the results in Section 8 that the space
of all bundle gerbe connections is an affine space for the vector
space $\Omega^1(Y)/ \pi^*(\Omega^1(M))$ of all one-forms on $Y$ modulo
one-forms pulled back from $M$.

\bigskip
\noindent{\bf   7. The curvature of a bundle gerbe connection. \hfill}

A bundle gerbe connection $\nabla$
is a connection so it has a curvature $F_\nabla$
which is a two-form on $Y^{[2]}$.
In the case that this is a trivial gerbe and the
connection is the tensor product connection then the curvature
can be written as
$$
 F_\nabla = \pi_2^* f - \pi_1^* f
\eqno(7.1)
$$
where $f$ is the curvature of the connection on the bundle over $Y$ and
the $\pi_i$ are the two projections $\pi_i \colon Y^{[2]} \to Y$.
We shall see in the next section that it is always possible to
find  an $f$ satisfying  equation (7.1). This is certainly true for a
bundle gerbe connection constructed
by a partition of unity argument as in Section 6.
The choice of such an $f$
we will call a curvature for the gerbe connection.

Consider now $df$. In the case that this is a trivial gerbe we have,
of course, that $df= 0$. More generally we have $dF_\nabla = 0$ so that
$$
\pi_1^* df = \pi_2^* df.
\eqno(7.2)
$$
I claim that this means that $df = \pi^*(\omega)$ for
some three-form $\omega$ on $M$.
To see this note that  a point of $Y^{[2]}$
is a pair $(x, z)$ where $\pi(x) = \pi(z)$
and a tangent vector to $Y^{[2]}$ at $(x, z)$
is a pair $(X, Z)$ with $X \in T_xY$ and
 $Z \in T_zY$  and $\pi_*(X) = \pi_*(Z)$.
Then equation (7.2) says that
$$
df(x)(X_1, X_2, X_3) = df(z)(Z_1, Z_2, Z_3)
\eqno(7.3)
$$
whenever $\pi_*(X_i) = \pi_*(Z_i)$ for $i = 1, 2, 3$.
Hence if $m \in M$ and $\xi_i \in T_mM$ choose
$x \in Y$ and $X_i \in T_xY$ such that $\pi(x) = m$
and $\pi_*(X_i) = \xi_i$ and define
$$
\omega(m)(\xi_1, \xi_2, \xi_3) = df(x)(X_1, X_2, X_3).
\eqno(7.4)
$$
Equation (7.3) shows that this definition is independent of the choice
of $x$ and the $X_i$.
Clearly $\pi^*(\omega) = df$ and moreover $\omega$ is closed.
We call $\omega/(2\pi i)$ the Dixmier-Douady form  of
the pair $(\nabla, f)$.

It will follow from the results in Section 8 that the
various choices in this construction do not
change the cohomology class of the
Dixmier-Douady form. We shall see in Section 11 that the
de Rham cohomology class of the Dixmier-Douady form is the
image in real cohomology of the
Dixmier-Douady class defined in Section 5. This proves
in particular that the Dixmier-Douady form is integral a
fact that also follows from the discussion of holonomy in
Section 12.

\bigskip
\noindent{\bf  8. A complex with no cohomology\hfill}

For a fibration $Y$ let $Y^{[p]}$ denote the $p$th fibered
product. There are projection maps $\pi_i \colon Y^{[p]} \to
Y^{[p-1]}$ which omit the $i$th element for each $i = 1 \dots
p$.  These define a map
$$
\delta \colon \Omega^q(Y^{[p-1]}) \to \Omega^q(Y^{[p]})
\eqno(8.1)
$$
by
$$
\delta(\omega) = \sum_{i=1}^p (-1)^i \pi_i^*(\omega).
\eqno(8.2)
$$
Clearly $\delta^2 = 0$ so that $\Omega^q(Y^{[*]})$
is a complex. We wish to show that this complex
has no cohomology.  This will then settle the question of
the existence of $f$  in Section 8 as we
have $F_\nabla \in \Omega^2(Y^{[2]})$
with $\delta(F_\nabla) = 0$ and hence
$F_\nabla = \delta(f)$ for some $f \in
\Omega^2(Y^{[1]}) = \Omega^2(Y)$.

Consider first the case that the fibration is trivial, say
$Y = M \times F$.  The general case will follow by a partition of unity
argument.  In this case $Y^{[p]} = M \times F^p$. Because the
notation is cumbersome at this point it will be
convenient to denote a  collection of $q$ vectors
$(X^1, \dots , X^q)$ just by $X$ and the action of a $q$ form
$\tau$ on these vectors by $\tau(X)$ rather than
 $\tau(X^1, \dots, X^q)$.
When we are dealing with vectors tangent to $F^{p+1}$ at
$f = (f_1, \dots, f_{p+1})$ then
each of the $X$ is  a collection of vectors
$(X_1, \dots X_{p+1}) $ where each $X_j $ is a $q$-tuple of
vectors in $T_{f_j}(F)$.
So, with these notational conventions we have
for $\omega \in \Omega^q(Y^{[p]})$
$$
\displaylines{
\delta(\omega)(m,f) ( \xi, ( X_1, \dots, X_{p+1}) )\hfill \cr
\quad\quad= \sum_i (-1)^i \omega(m, f_1, \dots,
\hat f_i, \dots, f_p)((\xi, ( X_1,
\dots, \hat{X_i}, \dots, {X_{p+1}})) )
\cr}
$$
$$
\eqno(8.3)
$$
where $\xi$ is a $q$-tuple of vector tangent to $M$ at $m$.
Now fix a point $f$ in $F$ and a $q$ tuple of vectors
$X$ in the tangent space at $f$ and define
$\rho \in \Omega^q(Y^{[p-1]})$  by
$$
\rho(m, f_1, \dots, f_p)(\xi, X_1, \dots, X_p)
 = \omega(m, f_1, \dots, f_p, f) (\xi, X_1, \dots, X_p, X).
\eqno(8.4)
$$
It follows from the fact that $\delta(\omega) = 0$ that
$\delta(\rho) = (-1)^{p+1}\omega$. This proves the
required result in the case
that $Y$ is the trivial fibration.   The general case is
proved by choosing an open cover $U_\alpha$ such that $Y$ is trivial
over each $U_\alpha$. Then let $\psi_\alpha$ be a partition of unity
subordinate to that cover.  Let $Y_\alpha$ by the
part of $Y$ sitting over $U_\alpha$ and similarly for $(Y^{[p]})_\alpha$.
Note that $(Y^{[p]})_\alpha = (Y_\alpha)^{[p]}.$
There are projection maps for each  $Y^{[p]} \to M$ and we can
pull the partition of unity back to any of these spaces. We will denote
it by the same symbol.  If we start with $\omega$ in $\Omega^q(Y^{[p]})$
then $\omega_{\vert U_\alpha} = \delta(\rho_\alpha)$
for some $\rho_\alpha$.
Hence we have
$$
\omega = \sum_\alpha \psi_\alpha \delta(\rho_\alpha) =
\delta (\sum_\alpha \psi_\alpha \rho_\alpha) = \delta(\rho)
\eqno(8.5)
$$
where $\rho = \sum_\alpha \psi_\alpha \rho$.

At $ p = 1$ we define $Y^{[0]} = M$ and
let $\delta \colon \Omega^q(M) \to \Omega^q(Y)$
be pull-back under $\pi$.
Exactness follows exactly as for the proof, in Section 7,
that there exists an $\omega$ such that
$\pi^*(\omega) = df$.

We can now confirm two facts stated earlier. The first
is the affine space structure of the space of all bundle
gerbe connections. If $\nabla$ and $\nabla'$ are
two bundle gerbe connections they clearly differ by a one-form
$\eta$ on $Y^{[2]}$ with $\delta(\eta) = 0$. So $\eta = \delta(\mu)$.
On the other hand $\nabla + \delta(\mu)$ is a bundle gerbe connection
for any $\mu$ on $Y$ and $\delta(\mu) = 0$ for such a $\mu$
precisely when $\mu$ is pulled back from $M$. This gives the
required result. The second fact is
the independence of the class of the Dixmier-Douady form
from various choices. The first is the choice of $f$
satisfying $\delta(f) = F_\nabla$. If $f'$ is another
such then $\delta(f - f') = 0$ and hence $ f - f' =   \pi^*(\rho)$
so that $df - df' = \pi^*(rho)$ and $\omega - \omega' = d\rho$.
The other choice is the choice of bundle gerbe connection.
If  we have  $\nabla' = \nabla + \delta(\mu)$ then we can choose
$f'$ so that $f = f' + d\mu$ and hence $df = df'$ so that $\omega = \omega'$.

\bigskip
\noindent{\bf  9. Gerbes, connective structures and curvings \hfill}

The relationship with the theory of gerbes discussed in [1]  is as follows.
For any open set $U \subset M$ let $C(U)$ be the
set of all sections of $Y$ which we want to think of as the
objects in a category, in fact in  a groupoid. If $s$ and $t$
are two such sections they define a section $(s,t)$ of $Y^{[2]}$ over
$U$ by $m \mapsto (s(m), t(m))$. The morphisms from
$s$ to $t$ we define to be the sections of the bundle $(s,t)^{-1}P$
over $U$. The composition is constructed from the
composition on $P$.  This construction defines a pre-sheaf
of groupoids. The sheafification of this presheaf gives rise to
a gerbe in the sense of  Brylinski [1].

In [1] Brylinski introduces the notion of connective structure
and curving. We indicate here how these are related to the bundle
gerbe connection and its curvature.
Let $\nabla$ be a bundle gerbe connection for the bundle gerbe $P$
over $Y \to M$. Let $U$ be an open subset of $M$ over which
$Y$ admits a section $s \colon U \to Y$. Denote by $\hat s$ the induced
map $Y_{|U} \to (Y_{|U})^{[2]}$ defined by $y \mapsto (y, s(\pi(y))$.
Then we have an isomorphism $P \simeq \hat s^{-1}P \otimes \hat s^{-1}P^*$
 defined by the product
$$
P_{(p,q)} \to P_{(s(\pi(p), p)}^* \otimes P_{(s(\pi(p)), q)}
\eqno(9.1)
$$
which trivialises $P$ over $Y_{|U} \to U$.
Consider the set $Co(s)$ of all connections $A$ on the bundle $\hat s^{-1}P$
such that $\delta(A) = \nabla$. This space of connections is an affine
space  for $\Omega^1_U$ the space of all $1$-forms on $U$ and hence
$Co(s)$ defines a  $\Omega^1_U$ torsor.
This torsor is a connective structure
in the sense of Brylinski.

Assume now that we have chosen a two-form $f$ on $Y$ such that
$\delta(f) = F_\nabla$ where $F_\nabla$ is the curvature of $\nabla$.
Then to any $A$ in $Co(P)$ we can define
a two-form $K(A)$ on $U$ by
$$
\pi^*(K(A)) = F_A - f
\eqno(9.2)
$$
where $F_A$ is the curvature of $A$. This equation makes sense because
$$
\delta(F_A - f) = F_{\delta(A)} - \delta(f) = F_\nabla - F_\nabla = 0.
\eqno(9.3)
$$

Finally notice that $\pi^*(dK(A)) = -df $
so that $dK(A) = - \omega$ so that,
up to sign, this is this is the curvature of the bundle gerbe.

The definition of morphism of bundle gerbes on $M$ in section 3
naturally gives rise to a notion of isomorphism.  A fundamental result
about gerbes is the theorem   that  the Dixmier-Douady class gives an
exact correspondence between elements of $H^3(M, \Z)$ and equivalence
classes of gerbes [1].  This is not true for bundle gerbes on
$M$ and bundle gerbe isomorphism; there are bundle gerbes on $M$
which are not isomorphic but which have the same Dixmier-Douady
class and hence define equivalent gerbes. Indeed it is not hard to show
that if  $(\alpha,
\beta, 1_M)$ is a morphism of bundle gerbes $(P, Y, M)$ and $(Q, X,
M)$ then the Dixmier-Douady classes of $(P, Y, M)$ and $(Q, X, M)$ are
the same.  In this example $X$ and $Y$ can be quite different. For
example if $Y$ admits a global section we can take $X$ to be the image
of that section and $Q$ the restriction of $P$ to $X$. A bundle gerbe
where the fibers are points clearly has Dixmier-Douady class zero and
we have already seen that a bundle gerbe where the fibration has a
section also has Dixmier-Douady class zero. This dependence of  bundle
gerbes on the choice of a fibration is nicely eliminated by the gerbe
concept.

\bigskip
\noindent{\bf  10. The tautological bundle gerbe \hfill}

Let $\omega/(2\pi i)$ be a form representing a class in
$H^3(M, \Z)$ where $M$ is $2$ connected.  We shall show how to
construct a fibration of  groupoids with  $\omega$ as its curvature.
Recall that if $\Sigma$ is an  oriented two sphere in $M$
the Wess Zumino Witten action is an element of $\C^\times$
associated to $\Sigma$ by extending $\Sigma$ to a ball $B$ in $M$
and defining
$$
\wzw(\Sigma) = \exp( \int_B {\omega}).
\eqno(10.1)
$$
Similarly if $\Sigma$ and $\Sigma'$ are two disks in $M$ with  common
boundary denote by $\wzw(\Sigma, \Sigma')$ the  Wess Zumino Witten
action of the sphere  formed by their union if it is given
the orientation of the first disk.

Fix a base point for $M$ and let $Y \to M$ be the path-fibration.
Then $Y^{[2]}$ consists of all pairs of paths beginning
at the basepoint and with the same endpoints. Define the fibre
of $P$ at such a point by taking all pairs
consisting of a piecewise smoooth surface with these two
paths as boundary and a non-zero complex number and
defining an equivalence relation
$$
(\Sigma, z) \sim (\Sigma', z')
\eqno(10.1)
$$
if $z = \wzw(\Sigma, \Sigma')z' $. Denote equivalence classes by
square-brackets.
Then the set of all equivalence classes
 forms a principal $C^\times$ bundle over $Y$. We
need to show that it is a bundle gerbe by constructing a product.

The product map $P_{(x, y)}\otimes P_{(y, z)} \to P_{(x, z)}$
is defined by
$$
[\Sigma, z] \otimes [\Sigma', z'] \to [\Sigma \cup \Sigma', zz'].
\eqno(10.2)
$$
This makes sense because $\Sigma$ and $\Sigma'$ have  half of
each of their boundaries (the curve $y$ ) in common.

We now show  that this bundle gerbe has a bundle gerbe connection
whose Dixmier-Douady form is $\omega/(2\pi i)$. We could perform
calculations analogous
to those in Section 2 however it is simpler to actually use those
calculations as follows. Consider the evaluation map
$$
\ev \colon Y \times [0,1] \to M
\eqno(10.3)
$$
and use it to define a closed  two-form $f = \int_0^1 \ev^*(\omega)$.
Note that $f/(2\pi i)$ is integral. We can now repeat the constructions
in Section 2 but restrict them to $Y^{[2]} \circ Y^{[2]}
\subset Y^{[2]} \times Y^{[2]}$.  This defines the bundle gerbe
with connection $\nabla $ and curvature $F_\nabla$ satisfying
$F_\nabla = \pi^*_2(f) - \pi^*_1(f)$. It is now an easy
calculation to show that if $\pi \colon Y \to M$ then
$df = \pi^*(\omega)$ as required.

\bigskip
\noindent{\bf   11. Deligne cohomology\hfill}

The Deligne cohomology of $M$ that we are interested in
is the total cohomology of the $\log$-complex
$$
0 \to C^\infty(M)^\times \to \Omega^1(M)  \to \dots \to \Omega^p(M) \to 0.
\eqno(11.1)
$$
Here the first non-zero map is the exterior derivative of the log
or $f \mapsto df/f$. If $p = 1$ then the elements of $H^1$ of this total
cohomology   are represented in Cech cohomology with
respect to an open cover by  pairs $(A_\a, \sigma_{\a \b})$
subject to the condition that
$$
A_\a - A_\b = \sigma_{\a\b}^{-1} d \sigma_{\a\b}.
\eqno(11.2)
$$
It is not hard to show that the elements of this cohomology
are equivalence classes of $\Cstar$  bundles with connection.

We shall show that in the case  $p=2$ that we can manufacture
a class in this total cohomology from a gerbe with connection and
curvature. A class in this cohomology will be a triple
$$
(f_\a, A_{\a\b}, g_{\a\b\c}).
\eqno(11.3)
$$
These have to satisfy
$$
A_{\a\b} - A_{\b\c} + A_{\c\a} =  g_{\a\b\c}^{-1}dg_{\a\b\c}
\eqno(11.4)
$$
and
$$
f_\a - f_\b = dA_{\a\b}.
\eqno(11.5)
$$

We have already a candidate for $g_{\a\b\c}$. For the other two
we first let $\nabla_{\a\b}$ be the pull-back connection
$$
\nabla_{\a\b} = (s_\a, s_\b)^* \nabla
\eqno(11.6)
$$
and then define
$$
A_{\a\b} = \sigma_{\a\b}^*(\nabla_{\a\b}).
\eqno(11.7)
$$
We also define
$$
f_\a = s_\a^* f.
\eqno(11.8)
$$
The first relation follows from the fact that a bundle gerbe
connection preserves the product. So
$$
\nabla_{\a\b} \otimes \nabla_{\b\c} = \nabla_{\a\c}.
\eqno(11.9)
$$
But the pull-back of
$\nabla_{\a\b} \otimes \nabla_{\b\c}$ with
$\sigma_{\a\b}\otimes \sigma_{\b\c} $ is
$$
A_{\a\b} + A_{\b\c}.
\eqno(11.10)
$$
On the other hand this is also the pull-back of $A_{\a\c}$ with
$\sigma_{\a\c}g_{\a\b\c}$ or
$$
A_{\a\c} + g_{\a\b\c}^{-1}dg_{\a\b\c}
\eqno(11.11)
$$
and the result follows.

The second relation follows from the equation
$$
F = \pi^*_{1} f - \pi^*_{2} f
\eqno(11.12)
$$
by pulling-back both sides with $(s_\a, s_\b)$.
On the LHS we get $dA_{\a\b}$ and on the RHS
$f_\a - f_\b$.

It now follows readily that $df_\a$ is just the restriction
of $\omega$ to $U_\a$ and moreover by standard double complex
arguments it also follows that the class defined by $\omega/( 2 \pi i)$
is the same as the class defined by $g_{\a\b\c}$.
So the class defined by the Dimxier-Douady  form  is the image in
$ H^3(M, \R)$ of the Dixmier-Douady class in $ H^3(M, \Z)$.

\bigskip
\noindent{\bf 12. Holonomy of a bundle gerbe connection over
a two-sphere\hfill}

If we calculate the Deligne cohomology for a
bundle gerbe with connection and
curvature whose base-manifold is a two-sphere we can show that it is
$U(1)$. The resulting number is a generalisation of holonomy. A simple
way of understanding what this is is to first consider the case that
the base manifold is not a two-sphere but in fact has $\pi_2(M) = 0$.
In that case consider an embedded two-sphere $\Sigma $ in
$M$ which is the boundary of some ball $B$ in $\Sigma$.
If $\omega/(2\pi i)$ is the Dixmier-Douady class then
$$
\wzw(\Sigma) = \exp( \int_B {\omega})
\eqno(12.1)
$$
is the holonomy of the bundle gerbe connection over $\Sigma$.
In general this is not a satisfactory solution as it analogous  to it
is like defining holonomy for a connection on a $\Cstar$ bundle
by integrating curvature over spanning disks.  In general we want
an answer intrinsic to the two-sphere in question.

To define such the holonomy choose a point in the
two-sphere $\Sigma$ and think of it as a disk
with boundary identified.  It is possible to lift this disk to a
disk $D$ in $Y$ whose buondary lies entirely in one fibre.
  If we fix a point $y$ in this fibre we define a map of the
fibre into $Y^{[2]}$ by  $y' \mapsto (y', y)$.  We can then
calculate the holonomy of the connection on $P$ around the image
of the boundary of $D$ under this map. Call this
$\hol(\partial D, \nabla)$. In
addition we can integrate $f$ over $D$ and form
$$
\hol(\partial D, \nabla)^{-1} \exp( \int_{D} f).
\eqno(12.2)
$$

We need to show that this is independent of various choices.
Let us leave the base point fixed for a moment and consider a
second lift $D'$. Then we can define a disk $\tilde D$ inside $Y^{[2]}$
by
$$
\tilde D = \{ (x, y) \mid x \in D, y \in D'\}.
\eqno(12.3)
$$
{}From the elementary properties of holonomy we have
$$
1 = \hol(\partial \tilde D, \nabla)^{-1} \exp( \int_{\tilde D} F_\nabla)
\eqno(12.4)
$$
and from the equation satisfied by $f$, (7.1) we conclude that
$$
\hol(\partial D, \nabla)^{-1} \exp(\int_{D} f)
= \hol(\partial D', \nabla)^{-1} \exp( \int_{D'} f).
\eqno(12.5)
$$

 To show that the holonomy is independent of the base point
is a similar type of calculation but more involved. Let $D_1$ and
$D_2$ be two lifts of $\Sigma$ with different basepoints. Let
$D_3$ denote the map of the cylinder into $Y$ which covers $\Sigma$
with the two basepoints removed and at each basepoint is coincides
with either the boundary of $D_1$ of $D_2$. Now by considering
each pair of $D_i$ respectively we define subsets of $Y^{[2]}$
by
$$
D_{ij} = \{ (x, y) \mid x \in D_i , y \in  D_j\}.
\eqno(12.6)
$$
Notice that topologically the union of these three cylinders
is a cylinder with each end in the diagonal inside $Y^{[2]}$.
On the diagonal the connection $\nabla$ is flat so it follows
that the integral of curvature $F$ over the union of the
$D_{ij}$ is $2\pi i $ times an integer.  Expanding out this
integral as before gives the required result.

 It is straightforward to check that
if $B$ is a three ball in $M$ then we have
$$
\hol(\partial B, \nabla) = \exp( \int_B {\omega})
\eqno(12.7)
$$
We can use (12.7) to prove that the
Dixmier-Douady form is  integral.  Indeed if $X \subset M$
is any three dimensional submanifold of $M$ consider a family of
three balls $B_r$ inside $X$ shrinking to a point as $r \to 0$.
Then the integral of $\omega$ over $ X - B_r$
is the holonomy over the boundary of $B_r$ but as
$B_r$ shrinks to a point this has limit $1$ and hence the
exponential of the integral of $\omega$ over  all of $X$ is
$1$. So the Dixmier-Douady form, $\omega/{(2\pi i)}$,  is integral.

\bigskip
\noindent{\bf 13. Topology of bundle gerbes \hfill}

So far we have ignored any properties that the fibering $Y \to M$
must satisfy for there to be a non-trivial bundle gerbe $P \to Y^{[2]}$.
However the examples
we have considered all have $Y \to M$ having infinite dimensional
fibres and we shall show now that this is a necessary condition.
We use the result from [3] that if $Y \to M$ is does not have
infinite dimensional fibres then any smooth choice of closed
$p$-form on the fibres is the restriction of a closed $p$-form
 on $Y$. If this is true then consider the two-form
$f$ on $Y$. Its restriction to each fibre is closed and
hence by the  theorem in [3] there exists a closed two-form $\rho$ on
$Y$ such that $f - \rho $ is vertical. But $f-\rho$
and $d(f - \rho)$ are both vertical so we can find a two-form
 $\mu$ on $M$ such that $\pi^*(\mu) = f - \rho$.
Finally $\pi^*(\omega) = df = df - d\rho = \pi^*(d\mu)$
so that $\omega - d\mu$ and the bundle gerbe has trivial Dixmier-Douady
class.

\bigskip
\noindent{\bf  Acknowledgements   \hfill}

The financial support of the Australian Research Council and
useful conversations with
Alan Carey are gratefully acknowledged. My thanks also to Dan Asimov,
David L. Johnson, Geoffrey Mess, Richard Palais, Lorenzo Sadun,
Jan Stevens and Joe Wolf who replied to my query in the Usenet
newsgroup sci.math.research for information about bundles
on $S^3$.

\bigskip
\noindent{\bf   References   \hfill}

\item{[1]} Brylinksi, J-L.:
{\it Loop spaces, Characteristic classes and Geometric Quantization.}
Birk\allowbreak h\"auser, Berlin, 1992.

\item{[2]} Mac Lane, S.:
{\it Categories for the working mathematician.}
Springer-Verlag, New York. 1971.

\item{[3]} Gotay, M.J.; Lashof, R.; \'Sniatycki, J. and Weinstein, A.;
Closed forms on symplectic fibre bundles. {\it Comment. Math. Helvitici.}
{\bf 58},
(1983),
617--621.

\bye